\documentclass[12pt]{amsart}

\usepackage{epsfig}
\usepackage{amsthm}
\usepackage{amssymb}
\usepackage{amsmath}
\usepackage{amscd}

\newtheorem {Theorem}   {Theorem}

\numberwithin{Theorem}{section}
\newtheorem {Corollary}[Theorem]{Corollary}
\newtheorem {Lemma}[Theorem]    {Lemma}
\newtheorem {Proposition}[Theorem]{Proposition}
\theoremstyle{definition}
\newtheorem{Definition}[Theorem]{Definition}

\theoremstyle{remark}
\newtheorem{Remark}[Theorem]{Remark}

\newtheorem{Claim}[Theorem]{Claim}

\expandafter\chardef\csname pre amssym.def at\endcsname=\the\catcode`\@
\catcode`\@=11
\def\undefine#1{\let#1\undefined}
\def\newsymbol#1#2#3#4#5{\let\next@\relax
 \ifnum#2=\@ne\let\next@\msafam@\else
 \ifnum#2=\tw@\let\next@\msbfam@\fi\fi
 \mathchardef#1="#3\next@#4#5}
\def\mathhexbox@#1#2#3{\relax
 \ifmmode\mathpalette{}{\m@th\mathchar"#1#2#3}%
 \else\leavevmode\hbox{$\m@th\mathchar"#1#2#3$}\fi}
\def\hexnumber@#1{\ifcase#1 0\or 1\or 2\or 3\or 4\or 5\or 6\or 7\or 8\or
 9\or A\or B\or C\or D\or E\or F\fi}

\font\teneufm=eufm10
\font\seveneufm=eufm7
\font\fiveeufm=eufm5
\newfam\eufmfam
\textfont\eufmfam=\teneufm
\scriptfont\eufmfam=\seveneufm
\scriptscriptfont\eufmfam=\fiveeufm

\catcode`\@=\csname pre amssym.def at\endcsname

\def    \eps    {\epsilon}

\newcommand{\R}{{\bf R}}

\newcommand{\Z}{{\bf Z}}

\newcommand{\p}{{\partial}}
\newcommand{\al}{{\alpha}}

\newcommand{\Om}{{\Omega}}
\newcommand{\om}{{\omega}}

\newcommand{\ga}{{\gamma}}

\newcommand{\la}{{\lambda}}
\newcommand{\si}{{\sigma}}
\newcommand{\Crit}{{\rm Crit}}

\newcommand{\Aa}{{\mathcal A}}

\newcommand{\Dd}{{\mathcal D}}

\newcommand{\Ll}{{\mathcal L}}
\newcommand{\Jj}{{\mathcal J}}

\newcommand{\Ff}{{\mathcal F}}

\newcommand{\Hh}{{\mathcal H}}

\newcommand{\Mm}{{\mathcal M}}

\newcommand{\Uu}{{\mathcal U}}

\newcommand{\Si}{{\Sigma}}

\newcommand{\Diff}{{\rm Diff}}

\newcommand{\NI}{{\noindent}}

\def    \I      {{\mathbb I}}

\def    \R  {{\mathbb R}}
\def    \Z  {{\mathbb Z}}

\def    \p  {\partial}

\begin{document}






\title[Length minimizing Hamiltonian paths]
{Length minimizing Hamiltonian paths for symplectically aspherical
manifolds}

\author[Ely Kerman]{Ely Kerman}
\author[Fran\c{c}ois Lalonde]{Fran\c{c}ois Lalonde}

\address{ Department of mathematics, University of Toronto ; 
D\'{e}partement de
math\'{e}matiques et de statistique, Universit\'{e} de Montr\'{e}al.}
\email{ekerman@fields.utoronto.ca ; lalonde@dms.umontreal.ca}

\thanks{The first author is supported by NSERC fellowship PDF-230728.
The second author is partially supported by a CRC, NSERC grant OGP 0092913
and FCAR grant ER-1199.}


\bigskip

\begin{abstract}
In this note we consider the length minimizing properties of
Hamiltonian paths generated by quasi-autonomous Hamiltonians on
symplectically aspherical manifolds. Motivated by the work of 
Polterovich and Schwarz  \cite{po,sc}, we study the role, in the 
Floer complex of the generating
Hamiltonian, of the global extrema which remain
fixed as the time varies. Our main result determines a  natural
condition  which implies that the corresponding path minimizes the
positive Hofer length. We use this to prove that a quasi-autonomous
Hamiltonian generates a length minimizing path if it has under-twisted 
fixed global extrema $P, Q$ and no contractible periodic orbits {\em with
period one} and action  outside the interval $[\mathcal{A}(Q),
\mathcal{A}(P)]$. This, in turn, allows us to
produce new  examples of  autonomous Hamiltonian flows which are length
minimizing for  all times. These constructions are based on the geometry
of coisotropic  submanifolds. Finally, we give a new proof of the fact
that quasi-autonomous  Hamiltonians generate length
minimizing paths over short time intervals.

\medskip
Nous \'etudions dans cette note les chemins de diff\'eomorphismes
engendr\'es par des hamiltoniens quasi-autonomes sur des vari\'et\'es
symplectiquement asph\'eriques. Motiv\'es par le travail de Polterovich  et Schwarz
\cite{po,sc}, nous examinons le r\^ole des extrema globaux et fixes au cours
du temps dans le complexe de Floer de l'hamiltonien. Notre principal
r\'esultat donne une condition suffisante naturelle pour que l'isotopie
hamiltonienne minimise la partie positive de la norme de Hofer. On en
d\'eduit qu'un hamiltonien quasi-autonome engendre une isotopie minimisant
la norme de Hofer s'il a des extrema $P,Q$ globaux fix\'es qui sont
``sous-tendus'' et n'a aucune orbite contractile { \em de p\'eriode}  $1$ et
d'action hors de l'intervalle
$[\mathcal{A}(Q),
\mathcal{A}(P)]$. Ceci nous permet de construire de nouveaux exemples
d'hamiltoniens autonomes qui induisent des flots minimisant pour tous
les temps. Ces constructions sont bas\'ees sur la g\'eom\'etrie des 
vari\'et\'es co-isotropes. On donne enfin une nouvelle preuve du fait que
tout Hamiltonien quasi-autonome engendre une isotopie minimisante sur un
intervalle suffisamment petit.  

\end{abstract}

\maketitle

\section{Introduction}\label{sec:intro}

On a closed symplectic manifold $(M, \omega)$, each time-dependent 
function
$H \in \mathcal{H}= C^{\infty}([0,1] \times M,\R)$ defines a 
time-dependent
Hamiltonian vector field $X_H$ via the equation
$$
i_{X_H} \omega = -dH.
$$
The
corresponding flow is denoted by $\phi^{t \in [0,1]}_H$ and  the space $Ham(M,\omega)$
consists of all the time-1 maps,
$\phi^1_H$, obtained in this manner.

Every path $h_t \colon [0,1] \to Ham(M,\omega)$ has a family of 
Hamiltonians $H \in \Hh$
which satisfy $h_t = \phi^t_H \circ h_0$. In \cite{ho1}, Hofer used these 
functions to define the 
length of the path $h_t$ by   
\begin{eqnarray*}
Length(h_t) &=& \int_0^1 \max_M H(x,t) - \min_M H(x,t) \,\,dt \\
            &=& \|H\|^+ + \|H\|^-, 
\end{eqnarray*}
which is independent of the choice of $H$. 

There is a unique Hamiltonian for each path that satisfies the 
normalization condition
\begin{equation*}\label{normal}
\int_M H(x,t) \om^n  = 0, \,\,\,\, \forall t \in [0,1].
\end{equation*}
For the normalized Hamiltonian $H$ of the path $h_t$, both $ \|H\|^+$ and 
$ \|H\|^-$ are non-negative and 
 provide different measures of length of $h_t$ called the positive and 
negative Hofer lengths, respectively.     

\begin{Remark} We will use both normalized and unnormalized Hamiltonians. 
Therefore, we will 
mention explicitly when this condition is assumed.
\end{Remark} 

A path is said to {\em minimize the (positive,negative) Hofer length} if 
there is no 
path in $Ham(M,\omega)$ with the same end points that is shorter (in the 
appropriate sense). 

Let us recall the properties which are necessary for $H$ to generate a 
path which minimizes the Hofer length.

\begin{Definition} A function $H \in \mathcal{H}$ is said to be
{\em quasi-autonomous} if it has at least one fixed global maximum $P \in 
M$ and one
fixed global minimum $Q \in M .$ In other words,
$$
H(t,P) \geq H(t,x) \geq H(t,Q) , \,\,\,\,\, \forall t \in [0,1]
\text{  and  } x \in M.
$$
\end{Definition}
Such fixed global extrema are clearly fixed points of the flow $\phi^{t \in [0,1]}_H$.

\begin{Definition}
A fixed point $x \in M$ of the flow $\phi^{t \in [0,1]}_H $ is said to be
{\em under-twisted} if, given any value $T
\in [0,1]$,  the linearized flow $D\phi^{t \in [0,T]}_H(x) \colon T_xM
\to T_xM$ has no non-constant T-periodic orbit. \footnote{We use the terminology
$T$-periodic orbit to mean any orbit that comes back to its initial position 
at
time $T$, i.e a closed orbit of period $T$; it should not be understood as {\em periodic in time}.} 
It is {\em generically under-twisted} if the origin is the 
only fixed
point of the flow $D\phi^{t \in (0,1]}_H$.  
\end{Definition}

\begin{Theorem}(\cite{bp,lm1,us})\label{nec}
If a Hamiltonian generates a length minimizing path then it must
be quasi-autonomous. Moreover, if there are
finitely many fixed global extrema, then at least one global
maximum and one global minimum must be under-twisted.
\end{Theorem}
In this paper we consider the converse question : When (for how long) 
is a  path length minimizing if it is generated by a
quasi-autonomous Hamiltonian with at least one (generically) under-twisted 
fixed
global maximum and minimum?  As described below, our approach relies 
on Polterovich's natural idea \cite{po} of applying Floer theoretic method to 
Hofer's geometry (see also Theorem~5.11 and Proposition~5.12 in Schwarz  \cite{sc}).
As in \cite{sc} and \cite{po}, 
we restrict ourselves here to the case when $(M,\omega)$ is
symplectically aspherical, i.e., $\omega |_{\pi_2(M)}=0=c_1|_{\pi_2(M)}.$ 
For such manifolds the action functional $\Aa_H \colon \Ll M \to \R$ is 
well-defined
on the space $\Ll M$ of contractible loops in $M$. It is given by

\begin{equation}\label{action}
\mathcal{A}_H(x(t)) = \int_0^1 H(t,x(t)) \, dt - \int_{D^2}
\bar{x}^* \omega,
\end{equation}
where $ \bar{x} \colon D^2 \to M$ satisfies $ \bar{x}|_{\partial D^2} = 
x$. 

Our first result is the following. 
 
\begin{Theorem}\label{main1}
Let $H \in \mathcal{H}$ be quasi-autonomous with a 
generically under-twisted fixed global maximum at $P$ and a 
generically under-twisted fixed global minimum at $Q$. 
If there are no  nonconstant contractible 1-periodic orbits of
$\phi^{t \in [0,1]}_H$ with action outside the interval
$[\mathcal{A}_H(Q),\mathcal{A}_H(P)]$, then $\phi^t_H$ is length 
minimizing over the time interval $[0,1]$. 
\end{Theorem}

In fact, Theorem~\ref{main1} will be shown to be a consequence of the 
following more general result.  
\begin{Definition}\label{def:dominate}
We say that $H \in \mathcal{H}$ {\em dominates} $K \in \mathcal{H}$   {\em 
at} $P$ 
if $H(t,x) \ge K(t,x)$ for all $t,x$
with equality at $x=P$ for all $t$. 
\end{Definition}
\begin{Theorem}\label{main-general}
Let $H \in \mathcal{H}$ have an under-twisted fixed global maximum at $P$. 
Assume that $H$ dominates some Hamiltonian $K$ at $P$, such that $P$ is 
a generically under-twisted fixed global maximum of $K$ and $\phi^{t \in [0,1]}_K$
has no contractible 1-periodic orbits with action greater than 
$\mathcal{A}_{K}(P) = \mathcal{A}_{H}(P)$. Then the path $\phi^{t \in [0,1]}_H$ generated by $H$ 
minimizes the
positive Hofer length.
\end{Theorem}
\begin{Remark}
A similar result also holds for the negative Hofer length of a Hamiltonian 
$H$ with
an under-twisted fixed global minimum at $Q$ (see Corollary \ref{cor:q}). 
\end{Remark}

   Note that our criterion depends only on the contractible 
periodic orbits {\em with period one}, not on periodic orbits with
intermediate periods. In the autonomous case, we get:

\begin{Corollary}\label{aut}
Let $H \colon M \to \R $ have a nondegenerate global maximum $P$ and 
minimum $Q$ which
are under-twisted. If there are no nonconstant contractible 1-periodic 
orbits of  $\phi^{t \in [0,1]}_H$ 
with action outside the interval $[H(Q), H(P)]$, then the path
$\phi^{t \in [0,1]}_H$ is length minimizing.
\end{Corollary}

   Of course, if there are no nonconstant contractible 1-periodic 
orbits of  $\phi^{t \in [0,1]}_H$, then it is  {\em a fortiori} length minimizing.

The papers \cite{lm2}, \cite{en},
\cite{mcsl} and
\cite{oh3} all include theorems  similar
to Corollary \ref{aut}, but which hold for more general symplectic 
manifolds.\footnote{Similar results were first obtained, for $\R^{2n}$, by 
Hofer in \cite{ho2} and 
Bialy-Polterovich in \cite{bp}.} However, they all require that there
be no nonconstant contractible 
periodic orbits with {\em period less than or equal} to $1$.

\medskip

Theorem~\ref{main1} also applies to quasi-autonomous Hamiltonians and 
allows for 
the existence of nonconstant contractible periodic orbits. This is 
significant 
because the existence of such orbits with all periods is a generic 
property for quasi-autonomous Hamiltonians. More precisely, if $H$ 
is suitably generic, then the Arnold conjecture 
implies that there are at least $SB(M)$ contractible $T$-periodic orbits 
for
each $T \in [0,1]$. Here $SB(M)$ denotes the sum of Betti numbers of $M$.
Hence, even when H is quasi-autonomous, there are, for most
manifolds,  contractible $T$-periodic orbits of $\phi^T_H$ other than the 
fixed global maximum and minimum of $H$. To assume that these other orbits 
are also constant for all $T \in (0,1]$, is clearly quite restrictive.

As an application of Theorem \ref{main1}, we will construct new examples 
of autonomous Hamiltonian flows which are length minimizing for all times. 
These constructions are based on the rich geometry of coisotropic 
submanifolds.

\medskip

Theorem \ref{main-general} also yields a proof of the following
result for the case of symplectically aspherical manifolds.
\begin{Theorem}(\cite{mc})\label{main2}
Let $H \in \mathcal{H}$ be quasi-autonomous. Then $\phi^t_H$
is length minimizing over $[0, \epsilon]$ for sufficiently small 
$\epsilon$.
\end{Theorem}
In \cite{lm2}, Lalonde and McDuff prove this result for any symplectic 
manifold 
but for the weaker notion of being length minimizing in the homotopy class 
of paths  
with fixed endpoints. Recently, McDuff in \cite{mc} has extended this to 
work with no 
assumption on the homotopy class of paths. We include the proof of the 
more restrictive case
above as a Floer theoretic interpretation of this fact.

\medskip
 The proof of Theorem \ref{main-general} relies on methods from 
Floer theory which were introduced to the study of Hamiltonian paths by 
Schwarz in \cite{sc} and Polterovich in \cite{po} \footnote{This is closely related to a
suggestion made by
Viterbo in  Remark 1.5B of Bialy-Polterovich \cite{bp}. Later, Oh extended these 
ideas to the action functional in \cite{oh1, oh2}.}. In \cite{oh3}, 
Oh also applies these methods, in a more general context, to
study Hamiltonian paths. 

Let us recall the basic idea as described in \cite{po}(Chapter 13).
Given a Hamiltonian $H$ which is suitably generic and has an under-twisted
fixed global maximum at $P$, one considers the role
of $P$ in the Floer complex of $H$. If one can show that $P$ is 
``homologically 
essential'' for the Floer complex, then certain perturbed 
pseudo-holomorphic 
cylinders must exist. An estimate for the energy 
of these cylinders then implies that the path generated by $H$  
minimizes the positive Hofer length. 
This strategy was used in \cite{oh3,po,sc} to prove length minimizing 
properties of paths generated by autonomous and quasi-autonomous 
Hamiltonians whose 
flows had no nonconstant contractible periodic orbits. 

In this paper, we simply refine this argument, in the 
simplest case of symplectically aspherical manifolds, in order to 
determine a length-minimizing criteria for Hamiltonian paths  
which allows for the existence of nonconstant contractible closed 
orbits. No new  methods are introduced here. Formally, 
the only new ingredient is that the notion of an element in the 
Floer complex ``being homologically essential'' is generalized to 
the notion of it being {\em homologically essential with respect to 
a filtration}: see the Definition~\ref{def:key}. This generalization
 is needed only to pass from Theorem~\ref{main1} to  
Theorem~\ref{main-general}, i.e to prove length minimizing 
properties for all Hamiltonians $H$ that dominate some other 
Hamiltonian $K$ satisfying the conditions in 
Theorem~\ref{main1}. As pointed out by Oh, this notion is already 
present in his earlier article \cite{oh3}, as an important step in his 
Lemma~7.8 (Non-pushing down Lemma). However, to the knowledge of the authors, 
it does not seem to be used there for the same purpose.

\medskip
Finally, we note that the results presented here yield a
description of the changes that must occur in the Floer chain complex in 
order for a generic quasi-autonomous
Hamiltonian to stop generating a length minimizing path. Let $H \in
\mathcal{H}$ be such a Hamiltonian with an under-twisted fixed
global maximum at $P$ and minimum at $Q$. Assume further that the 
contractible 
closed
orbits of $\phi^t_H$ are nondegenerate
for all $t \in (0,1]$. We then have the following generic picture
described by Floer in \cite{fl2}. There is a finite
set
$$
\{t_0=0,t_1,\dots,t_{k-1},t_k=1\} \subset [0,1]
$$
such that for $t \in (t_j ,t_{j+1})$ the number of contractible closed 
orbits 
of $\phi^t_H$ remains constant, as do the Conley-Zehnder
indices of the orbits. As $t$ increases through each $t_j$,
a pair of closed orbits of $\phi^{t_j}_H$, $\{x^-_j, x^+_j\}$, is
either created or destroyed. The Conley-Zehnder index of
$ x^+_j$ is one greater than that of $ x^-_j$ and the action of $x^+_j$
is greater than the action of $ x^-_j$ .

In this picture, the path $\phi^t_H$ is length minimizing at least
until the first $t_j$ at which a pair
 $\{x^-_j, x^+_j\}$ is created for which {\bf{either}} the Conley-Zehnder
index of $ x^+_j$ is one greater than that of $P$
{\bf{or}} the Conley-Zehnder index of $ x^-_j$ is one less than
that of $Q$. This can only happen at $t_j$ with $j>0$. Moreover, in the
first case it
is also necessary that the action of $ x^+_j$ be greater than that of
$P$. Similarly,
in the second case, the path is length minimizing until $ x^-_j$ has
action less than
$Q$.

\begin{Remark} \label{remark}  If the path
becomes non-minimal at some time, then it obviously remains non-minimal for all subsequent times.
Hence Theorem~\ref{main1} implies when a path
becomes non-minimal at some time $t'$,  such $t$-closed orbits must exist
for all times $t \ge t'$ as long as both global extrema remain
generically under-twisted. Since one is free to extend such a path by {\em any} Hamiltonian that keeps the minimum and maximum under-twisted, this is a surprising result. It is a consequence of the fact that 
our main theorem is stated in terms of  the non-existence of 
closed orbits of period one only.

\end{Remark}

What is described in this note extends to more general classes of symplectic 
manifolds  -- it will be developed further in the sequel \cite{kl2}. 
Geometric extensions of these ideas will also be studied in the paper
 \cite{la}.

\subsection{Organization of the paper}

In the next section, we recall the construction and properties of the 
Floer 
complex of a generic Hamiltonian. 
Following Schwarz, we describe in Section~\ref{sec:g-map},
how to identify the Floer complexes of normalized Hamiltonians which 
generate 
the same time-$1$ map. This identification is then used 
to prove the existence of a special cycle which represents the fundamental 
class in the 
Floer homology of a fixed normalized Hamiltonian. 
The notion of being 
``homologically essential with respect to a filtration'', which is equivalent to Oh's Non-pushing down Lemma, is given
in Section \ref{sec:hom-ess}. 
The proof of Theorem \ref{main-general} is contained in 
Section~\ref{sec:proof-general}. 
This result is then used in Section \ref{sec:proof-12} to prove Theorem 
\ref{main1} and Theorem \ref{main2}.
In Section~\ref{sec:infin}, we construct many new
examples of  autonomous functions which generate
Hamiltonian paths that are length minimizing for all time. Finally, we 
prove a technical claim (Claim 
\ref{claim:extra}) in the Appendix.

\bigskip

\section{Floer homology for symplectically aspherical symplectic 
manifolds}\label{sec:basics}

Let $\mathcal{H}_{reg} \subset C^{\infty}(S^1 \times M, \R)$ be the subset 
of Hamiltonians which are periodic in $t$ and whose contractible 
$1$-periodic
orbits are nondegenerate. In this section, we briefly recall the 
construction of the Floer chain
complex associated to each $G \in \mathcal{H}_{reg}$ when the symplectic 
manifold $(M,\omega)$ is symplectically aspherical. 
In particular, we focus on the algebraic aspects of these complexes that 
will be used later. The reader is referred to the sources 
\cite{fl1,fl2,fl3} for the full details.

\subsection{ The Floer chain complex }  A symplectic manifold  
$(M,\omega)$ is 
said to be symplectically aspherical if
$$
\omega |_{\pi_2(M)}=0 \,\,\, \text{    and   } \,\,\, c_1|_{\pi_2(M)}=0.
$$
Let $\Ll M$ be the space of contractible loops in $M$ 
 and consider the action functional $ \Aa_G \colon \Ll 
M \to \R$ defined by 
\begin{equation}\label{action}
\mathcal{A}_G(x(t)) = \int_0^1 G(t,x(t)) \, dt - \int_{D^2}
\bar{x}^* \omega,
\end{equation}
where $ \bar{x} \colon D^2 \to M$ satisfies $ \bar{x}|_{\partial D^2} = 
x$. The first condition to be symplectically aspherical implies  that the 
action functional
is well-defined.

The critical point set of $\Aa_G$ is equal to the (finite) set of 
contractible 1-periodic
orbits of $X_G$, i.e.,
$$\Crit(\Aa_G) =\{ x \in \mathcal{L}M \mid \dot{x} (t) = X_G(x(t)) \}.$$
The assumption that  $c_1|_{\pi_2(M)}=0$ implies that $\Crit(\Aa_G)$ is 
graded by the
Conley-Zehnder index, $\mu_{G}$. 

Here, we normalize $\mu_{G}$ as in \cite{se},
so that $\mu_{G}(x)$ is equal to the Morse index of $x$ whenever $x$ is a 
critical point of a $C^2$-small time-independent Hamiltonian. For an
under-twisted fixed maximum $P$ of $G$, $\mu_{G}(P)=2n$.  

The set $\Crit(\Aa_G)$ forms the basis of the Floer complex
$$
CF_*(G) =
\bigoplus_{k \in \Z} CF_k(G)
$$
where
$$
CF_k(G)=
\bigoplus_{x \in \Crit(\Aa_G),\, \mu_{G}(x)=k} \Z_2 x.
$$
We denote the obvious $\Z_2$-valued pairing for this
vector space by $\langle \,\,,\,\rangle$. The action of a 
chain $c \in CF(G)$ will be denoted by $\Aa_G(c)$ and defined as the 
maximum of the $\Aa_G$-values of the elements of $\Crit(\Aa_G)$ 
which appear in $c$ with a nonzero coefficient.

To construct the Floer boundary operator $\partial^G_{J_t}$
we first choose a family, $J_t =J_{t+1}$, of $\omega$-compatible
almost complex structures. For each pair $x,y
\in \Crit(\Aa_G)$, let $\Mm (x,y) = \Mm_{J_t}(x,y)$ be
the space of solutions 
$u \colon \R \times S^1 \to M$ of the equation
\begin{equation*}\label{dge}
\partial_s u + J_t(u)(\partial_t u -X_{G}(u)) =0
\end{equation*}
which satisfy the boundary conditions
\begin{equation*}\label{lim}
\lim_{s \to - \infty} u(s,t) = x  \,\,\text{    and  }
\lim_{s \to +\infty} u(s,t) = y,
\end{equation*}
and have finite energy
\begin{equation*}\label{en}
E(u)=\int_0^1 \int_{-\infty}^{+\infty} |\partial_s u|^2 \,ds\,dt.
\end{equation*}
If $J_t$ is suitably generic, then $\Mm_{J_t}(x,y)$ is a smooth
manifold of dimension $\mu_{G}(x)-\mu_{G}(y)$. We denote the set of 
these
generic $\omega$-compatible almost complex structures by
$\Jj_{\text{reg}}(G)$ and note that it is of second
category in the set of all $\omega$-compatible almost complex 
structures.

The boundary map $\partial^G_{J_t} \colon CF_k(G) \to CF_{k-1}(G)$ is now 
defined  by
$$
\partial^G_{J_t} x = \sum_{\mu_{G}(x) -\mu_{G}(y) =1} \tau(x,y)y,
$$
where $\tau(x,y)$ is the number $( mod(2))$ of elements in $\Mm (x,y)/\R$,
and $\R$ acts like $a \cdot u(s,t)= u(s+a,t)$.

For $J_t \in \mathcal{J}_{\text{reg}}(G)$, the boundary operator satisfies
$$
\partial^G_{J_t} \circ \partial^G_{J_t}=0.
$$
The homology of the complex $(CF_*(G), \partial^G_{J_t})$ 
does not depend on the choice of $J_t \in \mathcal{J}_{\text{reg}}(G)$ and
so we write it as $HF_*(G)$, the Floer homology of $G$.

\subsection{Homotopy chain maps and the action functional}
Let $G_s$ be a smooth homotopy $\R \to C^{\infty}([0,1] \times M, \R)$ 
such 
that
$$
G_s=\begin{cases}
 G_0 \text{ for } s<-\tau\\
 G_1 \text{ for } s>\tau,
\end{cases}
$$
for some $\tau >0$ and $G_0, G_1 \in \mathcal{H}_{reg}$.
Associated to each such homotopy is a homotopy chain map
\begin{equation*}
\sigma_{G_s} \colon CF_*(G_0) \to CF_*(G_1).
\end{equation*}

To construct this map we first choose a homotopy of families of
$\om$-compatible almost complex structures, $s \mapsto J_{s,t}$, such that
$$
J_{s,t}=
  \begin{cases}
    J^0_t \in \Jj_{\text{reg}}(G_0) & \text{ for }s < -\tau \\
    J^1_t \in \Jj_{\text{reg}}(G_1)& \text{ for }s> \tau.
  \end{cases}
$$
Then we consider the space $\Mm^s(x,y)=\Mm_{J_{s,t}}^s(x,y)$ of finite 
energy
solutions of the equation
\begin{equation*}\label{sdge}
\partial_s u + J_{s,t}(u) ( \partial_t u - X_{G_{s,t}}(u)) =0
\end{equation*}
which satisfy
\begin{equation*}\label{slim}
\lim_{s \to - \infty} u(s,t) = x\in \Crit(\Aa_{G_0}) \text{    and  }
\lim_{s \to +\infty} u(s,t) = y \in \Crit(\Aa_{G_1}).
\end{equation*}
For a suitably generic $J_{s,t}$, the space $\Mm_{J_{s,t}}^s(x,y)$ is a 
smooth 
manifold of dimension
$\mu_{G_1}(y) - \mu_{G_0}(x)$. The map $\sigma_{G_s}$ is then defined by
$$
 x \mapsto \sum_{ \mu_{G_0}(x) =
\mu_{G_1}(y)} \tau_s(x,y) y,
$$
where $\tau_s(x,y)$ is the number $(mod(2))$ of maps $u$ in the
zero-dimensional compact manifold $\mathcal{M}^s(x,y)$.

The next simple Lemma describes how the action changes under a homotopy 
chain map. 
\begin{Lemma}\label{a1}
Every $u \in \mathcal{M}^s(x,y)$ has energy
\begin{equation}\label{a1}
E(u) = \mathcal{A}_{G_0} (x) - \mathcal{A}_{G_1} (y) +
\int_0^1 \int _{-\infty}^{+\infty}
 \partial_s G(s,t,u(s,t))\,ds\,dt.
\end{equation}
\end{Lemma}
In the case of a linear homotopy, $G_s(x,t) = (1- \beta (s))G_0 +
\beta (s)G_1$, where $\beta$ is a smooth
nondecreasing function from $0$ to $1$, equation \eqref{a1} takes the
following simple form which will be used extensively later on.
\begin{equation}\label{a2}
E(u) = \mathcal{A}_{G_0} (x) - \mathcal{A}_{G_1} (y) +
\int_0^1 \int _{-\infty}^{+\infty}
 \dot{\beta}(s) (G_1 -G_0)(t,u(s,t))\,ds\,dt.
\end{equation}
Note that if $x \neq y$, then $E(u) >0$.

\subsection{ Homotopy chain maps and $HF_*(G)$ }

Every homotopy chain map induces an isomorphism in Floer homology and
so $HF_*(G)$ is independent of $G$. Moreover,
when $G$ is a $C^2$-small Morse function, the Floer 
complex
of $G$ (for a generic time-independent $J$) coincides with the Morse 
complex 
of $G$ (see \cite{hs}). Hence, for all $G \in \mathcal{H}_{reg}$,
the Floer homology $HF_*(G)$ is isomorphic to
the singular homology $H_*(M, \Z_2)$ via a map which preserves the 
grading.

\section{Comparing Floer complexes of normalized Hamiltonians 
generating the same time-$1$ map} \label{sec:g-map}

Given $\phi \in Ham(M, \omega)$, let
$\mathcal{F}(\phi)$ be the set of normalized Hamiltonians in 
$\Hh_{reg}$ whose time-$1$ map equals $\phi$.
For each pair of functions $F,G \in \mathcal{F}(\phi)$, there
is a natural chain isomorphism between the Floer complexes of $F$ and $G$.
These maps were studied by Seidel in \cite{se} for a very general class of 
symplectic 
manifolds and by Schwarz in \cite{sc} for the symplectically aspherical 
case. Let us
 recall the relevant details.

For $F,G \in \mathcal{F}(\phi)$, consider the loop
$g  = \phi^t_G \circ (\phi_F^t)^{-1}.$ This loop represents an element in 
$\pi_1(Ham(M,\omega), \I d)$ and acts on
$\Ll M$ by
$$
g(z(t))=g_t(z(t)).
$$
It is a straightforward consequence of the Arnold conjecture that the the orbit of $g$ on any point $x \in M$ is a contractible loop; thus $g$ takes $\Ll M$ back to itself.  Clearly, $g$ 
takes the 
$1$-periodic orbits of $X_F$ to
the $1$-periodic orbits of $X_G$. 
Since $\Crit(\Aa_F)$ and $\Crit(\Aa_G)$ are bases for the respective 
complexes, the induced map $g \colon CF(F) \to CF(G)$ is an isomorphism of 
$\Z_2$-vector spaces. 

If we set
$$
J_G = Dg_t^{-1} \circ J_F \circ Dg_t,
$$
then the moduli space
$\mathcal{M}_{J_F}^F(x,y)$ gets mapped bijectively, by $g$, to
$\mathcal{M}_{J_G}^G(g(x),g(y))$ for every
$x,y \in \Crit(\Aa_F)$, (see \cite{se} Lemma 4.3). Thus, $g$ is also a 
chain complex isomorphism
from  $(CF(F),\partial^F_{J_F})$ to $(CF(G),\partial^G_{J_G})$.

Moreover, for symplectically aspherical manifolds, the chain complex 
isomorphism $g$ preserves 
both the grading by the Conley-Zehnder index and the filtration by the 
action functional,(see \cite{sc} Theorem 1.1).

\medskip

Let $G \in \mathcal{H}_{reg}$ be a normalized Hamiltonian. The fact that 
the Floer 
complex of any other function $F \in \Ff (\phi^1_G)$ is identical to 
$(CF(G),\partial^G_{J_G})$ 
can be used to prove the following result.

\begin{Lemma} \label{le:comparison} For any $F \in \Ff
(\phi^1_G)$ there is a cycle $v \in CF_{2n}(G)$
 which generates $HF_{2n}(G)$ and satisfies $\mathcal{A}_{G}(v) \le 
\|F\|^+$.
\end{Lemma}

\begin{proof}
Let $f$ be a Morse function. For a sufficiently small 
$\epsilon >0$ and a suitably
generic time-independent almost complex structure, the Floer chain complex 
of $\epsilon f$ is well-defined
and equal to the Morse complex of $\epsilon f$. In particular, the cycle 
$\sum_i p^i_{max}$, given by the sum of the local 
maxima of $\epsilon f$, generates $HF_{2n}(\epsilon f)$.

The linear homotopy between $\epsilon f$ and $F$ yields the homotopy chain 
map
$$
\sigma \colon CF( \epsilon f ) \to
CF(F).
$$
Hence, the cycle $v^F = \sigma (\sum_i p_i^{max})$ represents the 
fundamental class in $HF_{2n}(F)$.

Let $v^F_i \in \Crit(\Aa_F)$ appear in $v^F$ with nonzero coefficient. 
Then, for some $j$, there must be an element
$ u \in \mathcal{M}^s( p_j^{max},v^F_i)$. Equation \eqref{a2} then yields 
the inequality
$$
 \mathcal{A}_F(v^F_i) \leq \mathcal{A}_{\epsilon f}( p_j^{max}) + \int_0^1 
\max_M \{F(x,t)-\epsilon f(x)\} \,dt.
$$
For sufficiently small $\epsilon > 0$ this implies that 
$\mathcal{A}_F(v^F_i) \leq \|F\|^+$. Therefore, 
$\Aa_F(v^F) \leq \|F\|^+.$

By the discussion above, the chain isomorphism from $(CF(F), 
\partial^F_{J_F})$ to $(CF(G), \partial^{G}_{J_G})$, 
induced by the loop $g=\phi^t_G \circ (\phi^t_F)^{-1}$, preserves indices 
and actions. Hence, the cycle $v=g(v^F)$ 
also generates $HF_{2n}(G)$ and satisfies $\Aa_G(v) \le \|F\|^+$.
\end{proof}

\section{ Homologically  essential with respect to a 
filtration}\label{sec:hom-ess}

Let $(C, \p )$ be a differential complex over $\Z_2$ which is freely 
generated by a basis $\mathcal B$.
\begin{Definition}\label{def:hom-ess}
Let $\ga \in H_*(C, \p)$ be a nontrivial homology class. An element
$b \in \mathcal B$ is said to be {\em homologically essential for the 
class $\ga$} if 
it appears with nonvanishing coefficient in every representative of $\ga$
\end{Definition}
In this section we refine this notion in the presence of the following 
additional structure.
\begin{Definition}  A {\em filtration} on $(C, \p)$
is a map
$\mathcal{V} : \mathcal B \to \R$ such that
$$
\mathcal{V} (b) > \mathcal{V} (b')
$$
whenever $b,b' \in \mathcal B$
satisfy $\langle b', \p b \rangle \neq 0$.  For an element $c \in C$, we 
define
$\mathcal V (c)$ to be the maximum of the values of $\mathcal V$ over the
elements of $\mathcal B$ that appear with non-vanishing coefficient in 
$c$.
\end{Definition}

A Morse function together with its Morse complex, and an action functional 
together with its Floer complex,
are obvious examples of filtrations.

\begin{Definition} \label{def:key} Let $(C, \p )$ be a complex with
filtration $\mathcal{V}$, and let $\ga \in H_*(C, \p)$ be a nontrivial 
homology class.
An element $b \in \mathcal B$ is said to be {\em homologically essential 
for $\ga$ with respect to
$\mathcal{V}$ } if the following two conditions are
satisfied:

\medskip
\NI
1) there is an element representing $\ga$ of the form $b + v$ such that
$\langle b, v \rangle = 0$ and $\mathcal V (v) < \mathcal V (b)$;

\medskip
\NI
2) whenever $\langle b, \p d \rangle \neq 0$, then $\mathcal V (\p
d - b) \ge \mathcal V (b)$.
\end{Definition}

To better understand the origin of these chain-level conditions, consider 
a 
Hamiltonian $H \in \Hh_{reg}$ with an under-twisted fixed global maximum 
at $P$. 
We will show (in Proposition \ref{prop:above}) that, under a generic
 nondegeneracy assumption, $P$ always satisfies the first condition to be 
homologically essential for the generator of $HF_{2n}(H)$ with respect to 
$\Aa_H$. 
This will follow from the fact that $H$ can be dominated at $P$ by a 
function $G^H$ for 
which $P$ satisfies the first condition to be essential for the generator 
of $HF_{2n}(G^H)$ with respect to $\Aa_{G^H}$. 

On the other hand, assume that $H$ dominates a function $G_H \in 
\Hh_{reg}$ for which 
$P$ is under-twisted and not in the image of the boundary operator 
$\partial^{G_H}$ 
(e.g.  $G_H$ has no $1$-periodic orbits with $\Aa_{G_H}$-value greater 
than $\Aa_{G_H}(P)$). 
We will also show that this implies that $P$ satisfies the second criteria 
of being homologically 
essential for $HF_{2n}(H)$ with respect to $\Aa_H$ (see Proposition 
\ref{prop:below}). 

Note that if $H$ satisfies the assumptions of $G_H$ itself, then $P$ is 
homologically  
essential for $HF_{2n}(H)$ in the sense of Definition \ref{def:hom-ess}.

\medskip

The next proposition is the main result of this section.

\begin{Proposition} \label{prop:ess-implies-min}
Let $G \in \mathcal{H}_{reg}$ be a normalized Hamiltonian with an 
under-twisted 
fixed global maximum at $P$. If $P$ is homologically essential for 
the unique generator of $HF_{2n}(G)$ 
with respect to $\Aa_{G}$, then the path generated by $G$ minimizes the 
positive Hofer
length, i.e., 
$$
\|G \|^+ \le \|F \|^+ \,\,\,\,\text{  for all   }\,\,F \in \Ff(\phi^1_G).
$$
\end{Proposition}

\begin{proof}
By our hypothesis on $P$, there is a  chain $w \in CF_{2n}(G)$ such that 
$[P + w]$ generates
$HF_{2n}(G)$ and $$\mathcal A_{G}(w) < \mathcal A_G (P).$$
If there is an $F \in \Ff(\phi^1_G)$ such that $\|F \|^+ < \|G \|^+$,  
then by Lemma~\ref{le:comparison} there is a
cycle $v$ that generates $HF_{2n}(G)$ and satisfies
$$
\mathcal{A}_{G}(v) \le \|F\|^+ < \|G \|^+ = \mathcal A_G(P).
$$
The two cycles $P+w$ and $v$ represent the same class, therefore
$P + (w - v)$ is exact. The second criteria for $P$ to be homologically 
essential for this class
with respect to $\mathcal A_G$ then implies that $\mathcal A_G (w-v)$ must 
be larger than 
or equal to $\mathcal A_G (P)$. This is a contradiction.
\end{proof}

\section{ Proof of Theorem~\ref{main-general}}\label{sec:proof-general}

In proving Theorem~\ref{main-general}, we may assume that the functions 
$H$ and $K$ 
have some additional properties.

First, we may assume that $H$ is normalized. This is clear, since the 
positive 
Hofer length of the path $\phi^t_H$ is independent of the normalization of 
$H$.
(When $H$ is normalized we just know that the positive Hofer length of 
$\phi^t_H$ equals 
$\|H\|^+$.)

We may also assume that $H$ and $K$ satisfy certain generic nondegeneracy 
assumptions.
Here is the precise statement.
\begin{Claim}\label{claim:extra}
It suffices to prove Theorem~\ref{main-general} for $H$ and $K$ with the 
following additional properties: 
\begin{enumerate}
\item $H,K \in \Hh_{reg},$ 
\item For all $t \in [0,1]$, $P$ is the unique global maximum of both 
      $H(t,x)$ and $K(t,x)$, and is nondegenerate (as a critical point). 
\end{enumerate} 
\end{Claim}
The proof of this claim is technical and more or less standard. It is 
therefore deferred 
to an appendix.

\medskip

Since we can assume that $H$ and $K$ are in $\Hh_{reg}$, we can consider 
their Floer 
complexes. By Proposition \ref{prop:ess-implies-min}, we will then be done 
if we can 
show that $P$ is homologically essential for the generator of $HF_{2n}(H)$ 
with respect 
to $\Aa_H$.  

\medskip

Let's begin with the first condition to be  homologically
essential.
\begin{Proposition}\label{prop:above}
Let $H \in \mathcal{H}_{reg}$ have a fixed under-twisted global maximum at
$P$ such that $P$ is the unique global maximum of $H(t,x)$ and is 
nondegenerate, 
for all $t \in [1,0]$. Then $P$ satisfies the first condition to be 
homologically 
essential for the unique generating class of $HF_{2n}(H)$ with respect to 
the filtration $\mathcal A_H$. 
That is to say:
$HF_{2n}(H)$ is generated by $[P+w]$ for some $w \in CF_{2n}(H)$ with
$\langle w, P\rangle=0$ and $\mathcal{A}_H (w)
< \mathcal{A}_H(P)$. 
\end{Proposition}
\begin{proof}
First we construct a function $G^H$ which dominates $H$ at $P$.
Let $f_P \colon M \to \R$ be an
autonomous Morse function for which $P$ is the only critical point with 
Morse
index $2n$ and $f_P(P) =0$. For
sufficiently small $\epsilon >0$ the function
$$
G^H(t,x) = H(t,P) + \epsilon
f_P(x)
$$
satisfies $ G^H(t,x) \geq
H(t,x),$ with equality only at $P$. This follows from the uniqueness and
nondegeneracy conditions on $P$.
Note that $P$ is also an under-twisted fixed global maximum of $G^H$ and 
$\mathcal{A}_{G^H}(P)=\mathcal{A}_{H}(P)$.

The Floer complex of $G^H$ is equal to the Floer complex of $\epsilon 
f_P(x)$
with the actions shifted upward by $
\int_0^1 H(t,P) \, dt$. For sufficiently small $\epsilon $, and a suitably
generic time-independent almost complex
structure, the Floer complex of $\epsilon f_P$ coincides with its Morse
complex. By our choice of the
function $f_P$, this means that $HF_{2n}(G^H)$ is generated by the class 
$[P]$.

The linear homotopy from $G^H$ to $H$ yields the homotopy chain map
$$
\sigma_{G^H} \colon CF_*(G^H) \to CF_*(H).
$$

Since this induces an isomorphism in homology, the class 
$[\sigma_{G^H}(P)]$ generates
$HF_{2n}(H)$.

Now, $P$ is an under-twisted fixed maximum of $H$, so $\mu_H(P) = 2n$.
Consider an element
$ u
\in
\mathcal{M}^s(P,P)
$. For each such
$u$, equation \eqref{a2} yields the inequality
$$
E(u) = \mathcal{A}_{G^H} (P) - \mathcal{A}_H (P) +
\int_0^1 \int _{-\infty}^{+\infty}
 \dot{\beta}(s) (H -G^H)(t,u(s,t))\,ds\,dt
$$
which implies that $E(u) \leq 0$. Thus, $ \mathcal{M}^s(P,P)$
consists of just the constant map, $u(s,t) \mapsto P$, and
$$
\sigma_{G^H}(P) = P + \sum_{i=1}^k m_{y_i} y_i=P+w.
$$
Here, $\{P, y_1, \dots ,y_k\}$ are the 1-periodic orbits of $H$
with $\mu_H=2n$.

To finish the proof we must show that any $y_j$ which appears in
$w$ with nonzero coefficient, must satisfy $
\mathcal{A}_H(y_j) < \mathcal{A}_H(P)$. For each element
$ u \in \mathcal{M}^s(P,y_j)$, equation \eqref{a2} yields
$$
0 < \mathcal{A}_{G^H}(P)- \mathcal{A}_H(y_j)+
\int_0^1 \int_{-\infty}^{\infty}
\dot{\beta}(s)(H-G^H)(t,u(s,t))\,ds\,dt.
$$
Hence, $\mathcal{A}_H(y_j)<\mathcal{A}_{G^H}(P)=\mathcal{A}_H(P)$.
\end{proof}

\medskip

Finally, we use the function $K$ to prove that $P$ satisfies the second 
condition 
to be homologically essential for the generator of $HF_{2n}(H)$ with 
respect to $\Aa_H$. 
\begin{Proposition}\label{prop:below}
If $\partial^{H} a = P+z$ for some $a \in 
CF_{2n+1}(H)$
and $z \in CF_{2n}(H)$, then
$\mathcal{A}_{H}(z) \geq \mathcal{A}_{H}(P)$.
\end{Proposition}

\begin{proof}
Let $\partial^{H} a = P+z$ be given.
Consider the linear
homotopy from $H$ to $K$. The induced
homotopy chain map $\sigma$ satisfies 
$$
\sigma (\partial^{H} a) =
\partial^{K} \sigma (a).
$$
Thus, we have 
$$
 \sigma(P) + \sigma(z) = \partial^{K} \sigma (a) .
$$
Since $K$ has no periodic orbits with action larger than $\Aa_K(P)$, $P$
occurs with coefficient zero in
$\partial^{K} \sigma (a).$ It is also straightforward to show that $P$
occurs in $\sigma(P)$ with coefficient $1$ (i.e., $\mathcal{M}^s(P,P)$
contains only the constant map). This means that $P$ must occur in
$\sigma(z)$ with coefficient $1$. In other words, if $z =\sum_i b_i z_i$,
then for some $z_i$ there exists  $u \in \mathcal{M}^s(z_i,P)$. Again by 
\eqref{a2}, this implies
that
$$
0< \mathcal{A}_{H} (z_i) - \mathcal{A}_{K} (P) +
\int_{-\infty}^{+\infty} \int_0^1
\dot{\beta}(s) (K-H)(t,u(s,t)) \,dt \, ds
$$
and   so
$$
\mathcal{A}_{H}(z) \ge
\mathcal{A}_{K}(P) = \mathcal{A}_{H}(P).
$$
\end{proof}

\section{ The Proofs of Theorem \ref{main1} and Theorem \ref{main2} 
}\label{sec:proof-12}

\subsection{Proof of Theorem \ref{main1}}
Consider a Hamiltonian $H \in \Hh$ with a fixed under-twisted global 
minimum at $Q$ and let
$$
\bar{H}(t,x)=-H(t,\phi^t_H(x)).
$$ 
This Hamiltonian generates the path $(\phi^t_H)^{-1}$ and has $Q$ as an 
under-twisted fixed global 
maximum. 
Note that the positive Hofer length of $(\phi^t_H)^{-1}$ is equal to the 
negative Hofer length of
$\phi^t_H$. Hence, applying Theorem \ref{main-general} to $\bar{H}$ 
yields:
\begin{Corollary}\label{cor:q}
If $\bar{H}$ dominates some Hamiltonian $L$ at $Q$ such that
$Q$ is a generically under-twisted fixed global maximum of $L$ and the 
flow of $L$
has no contractible 1-periodic orbits with action greater than 
$\mathcal{A}_{L}(Q) = -\mathcal{A}_{H}(Q)$, then the path generated by $H$ 
minimizes the
negative Hofer length.
\end{Corollary}
   
Theorem~\ref{main1} now follows immediately: simply take $H$ for $K$ in 
Theorem \ref{main-general}
and $\bar H$ for $L$ in Corollary \ref{cor:q}.

\subsection{Proof of Theorem \ref{main2}}

The flow of $H(t,x)$ is length minimizing for $t \in
[0,\epsilon]$ if and only if the flow of $\epsilon H( \epsilon t, x)$
is length minimizing for $t \in [0,1]$. Let $P$ and $Q$ be the 
fixed global maximum and minimum of $H$, respectively.
For all sufficiently small $\epsilon > 0$ these will be under-twisted
fixed global extrema of $\epsilon H( \epsilon t, x)$. 

Let $f_P \colon M \to \R$ be a Morse function for which $P$ is the
only critical point with Morse index $2n$ and $f_P(P)=0$. We also 
assume that $f_P$ is sufficiently $C^2$-small so that the only 
$1$-periodic 
orbits of $X_{f_P}$ are the critical points of $f_P$.

Set
$$
K=\epsilon H( \epsilon t, P) + f_P.
$$
Clearly, $P$ is a generically under-twisted fixed global maximum of $K$, 
and
$\epsilon H( \epsilon t, x)$ dominates $K$ at $P$ for sufficiently small 
$\epsilon$. 
Also, the only $1$-periodic orbits for the flow of $K$ are
the critical points of $f_P$, so there are also no $1$-periodic orbits 
with
action greater than $\Aa_K(P)$. Theorem \ref{main-general} then implies 
that the
path generated by $\epsilon H( \epsilon t, x)$ minimizes the positive 
Hofer 
length.

A similar argument applied to $\epsilon \bar{H}( \epsilon t, x)$ implies
that the path generated by $\epsilon H( \epsilon t, x)$ also minimizes the 
negative Hofer 
length. 
 
\section{Coisotropic manifolds and Hamiltonians generating length
minimizing  paths for all times}\label{sec:infin}
 Let $N$ be a closed manifold endowed with a closed 2-form $\om_N$ of
constant rank. Its kernel is then an integrable distribution that gives 
rise to
a foliation ${\mathcal F}$ on $N$. Up to local diffeomorphisms, there is a 
unique
way of considering $(N, \om_N)$ as a coisotropic submanifold. Indeed, let
$\rho: E
\to N$ be the vector bundle whose fiber at each point $q \in N$ is the 
dual $K^*_q$ of
the kernel $K_q$ of $\om_N$ at $q$. Choose a distribution $\Hh$ on $N$
transversal to $\Ff$, and define a form on each $T_{Z(q)}E, q \in N$,  by
$$
(d\rho)^*(\om_N) + \Pi^*(\om_{can})
$$
where $Z:N \to E$ is the zero section and the map $\Pi = (\Pi_1, \Pi_2):T_{Z(p)}E \to K \oplus K^*$ is the projection
induced by $\Hh$ (i.e $\Pi_2$ is the projection on the fiber of $E$ and $\Pi_1$ is the composition of the $d\rho$ with the projection induced by $\Hh$). This is a fibrewise symplectic form on the vector bundle 
$TE |_N$
whose restriction to $TN$ is $\om_N$. By the Moser-Weinstein theorem, 
there is
an extension of that form to a symplectic form $\om$ on some neighborhood 
$\Uu$
of $N$ in $E$, which is unique up to diffeomorphisms on smaller 
neighborhoods
whose 1-jet act as the identity on $TE |_N$. The submanifold $N$ is then 
coisotropic in $\Uu$.

  Let us now consider a special case of this construction. Let $W$ be a 
closed
manifold that admits a metric $g$ of non-positive
curvature, and let
$$
W \hookrightarrow  N \stackrel{\pi}{\to} B
$$
be a smooth fibration with structure group $\Diff_{g}(W)$, the group of
$g$-isometries of $W$. Assume also that the transition maps $\phi_{i,j}: 
V_i
\cap V_j \to \Diff_g(W)$ are locally constant. Thus each fiber $W_b$ is
equipped with a metric $g_b$ of non-positive curvature. As before,
$\om_N \in \Om^2(N)$ is any closed 2-form whose kernel at each point is 
the
tangent space to the $W$-fiber at that point. By differentiating the 
transition
functions of the bundle
$W \hookrightarrow  N \to B$, it can be extended to a $T^*W$-bundle
$$
T^*W \hookrightarrow  N'=E \stackrel{\pi'}{\to} B.
$$
There is a symplectic
form $d \la_b$ on each fiber $T^*(W_b)$ that extends to a closed form 
$\tau$ on
$E$. Since the transition functions are locally constant, one can choose
$\tau$ such that it coincides with $d \la$ in each local chart $V_i \times
T^*W$ (where $V_i$ is an open subset of $B$). The form
$\rho^*(\om_N) + \tau$, is then symplectic on $E$, where $\rho$ is now the 
projection $N' \to N$. 
This form also restricts to $\om_N$ on $N$.

Since the structure group preserves the metric, one may
define a Hamiltonian
$$
H: \Uu \to \R
$$
given by $H(b,p,q) = f( \| p \|^2) $ in local coordinates, where $f$ is 
the
identity map near $0$ and becomes constant for values larger than some 
$\eps >
0$. Extend $H$ outside $\Uu$ by the constant map. Since in local 
coordinates
$V_i \times T^*W$, the symplectic structure is $\rho^*(\om_N) + d\la$,
where $\rho^*(\om_N)$ has each $T^*W$-fiber as a kernel, the flow of $H$ 
is 
the geodesic flow along each leaf of the isotropic foliation of the 
coisotropic
submanifold
$N$. Because the metric has non-positive curvature, there is no 
non-constant
contractible closed orbit. Since this Hamiltonian can be approximated by 
an
autonomous Hamiltonian that is generically under-twisted at its fixed
maximum and minimum and has no non-trivial periodic orbit, we then have, 
as a
consequence of Theorem~\ref{main1} and the Lemma \ref{near} of Section~
\ref{sec:appendix}:

\begin{Corollary} Let $(M, \om)$ be an aspherical symplectic manifold that
contains a coisotropic submanifold of the above form. If the fundamental 
group
of each leaf injects in the fundamental group of $M$, the flow of $H$ is 
length
minimizing for all times.
\end{Corollary}

This result is very likely valid in non-aspherical manifolds as well,
though we will not present a general proof here. This corollary admits 
the
following two extreme cases.

\medskip
\noindent
{\bf Examples}. (1) The submanifold $N$ is a Lagrangian submanifold. In 
this
case, we recover the constructions due to Schwarz in \cite{sc}. It is
interesting to note that this can also be derived in the Lagrangian case
from previous results obtained by very different techniques that do
not involve Floer's homology. Actually, if $L \subset M$ is a Lagrangian
submanifold, the conditions that $L$ admits a metric with non-positive
curvature and that its fundamental group injects into the fundamental 
group of
$M$ means that the embedding of $L$ lifts to an embedding of the
contractible space $\tilde L$ into $\tilde M$. By Theorem 1.4.A in
Lalonde-Polterovich \cite{lp}, this implies that any neighborhood of 
$\tilde L$
has infinite capacity (there exist embedded balls of arbitrarily large
capacity) -- this was incidentally used in \cite{lp} to prove that no
Hamiltonian isotopy can disjoin $L$ from itself. Using now the same 
techniques
as in Lalonde-McDuff
\cite{lm2}, Lemma 5.7, one constructs an autonomous Hamiltonian flow
$\phi_t$ whose lift to the universal cover disjoins balls $B_t$ that have
capacity equal to the energy of $\phi_t$. This implies that $\phi_t$ is 
length
minimizing for all times by the energy-capacity inequality. Note that the
aspherical condition is not needed here.

\medskip
   (2) The submanifold $N$ is a hypersurface. This is a simple case since 
any
1-dimensional manifold admits a flat metric. Note that our construction of
coisotropic submanifolds implies both the stability of the coisotropic
submanifold and the non-existence of closed contractible orbits. This
is what we must require here, i.e the right statement is the following:
the existence of a length minimizing autonomous Hamiltonian flow is 
guaranteed
as soon as there exists in
$(M, \om)$ a stable hypersurface whose closed orbits are all 
non-contractible (by stable, we mean here that there is a neighbourhood of $N$ of the form 
$[-\eps, \eps) \times N$ in which each leaf $\{s\} \times N$ has a characteristic flow conjugated to the one on $N = \{0\} \times N$).
This is the case for instance for fibered symplectic manifolds over a surface of 
strictly
positive genus
$$
M \hookrightarrow (S, \Om) \to \Si
$$
i.e fibered manifolds equiopped with a symplectic form $\Om$ whose restriction to each
$M$-fiber is non-degenerate. The inverse image $N=\pi^{-1}(\ga)$
of any non-contractible loop in
$\Si$ yields a hypersurface whose characteristic foliation is transverse 
to the
$M$-fibers.  A closed orbit of such a foliation must necessarily project 
to a
non-vanishing multiple of $\ga$, and is therefore non-contractible. Note 
that
the symplectic structure in a neighborhood of the hypersurface has the 
form
$\Om |_N + d(t \al)$ where $t$ is the coordinate in the normal direction 
to $N$
and
$\al \in \Om^1(N)$ is any 1-form that does not vanish on the 
characteristic
directions. Choosing $\al$ as the pull-back of $d \theta \in \Om^1(\ga)$, 
the
symplectic form becomes $\Om |_N + dt \wedge d \theta$, for which $N$
is stable as required.

\medskip
   (3) One may construct lots of examples in between, i.e.
coisotropic submanifolds of dimensions $n < d < 2n-1$ that satisfy the
conditions of our construction. The following example shows that the 
existence
of length minimizing paths is stable under pull-backs via certain 
symplectic
fibrations. Let
$$
M \hookrightarrow (S, \Om) \stackrel{\pi}{\to} (B, \si)
$$
be a fibration with compatible symplectic forms.
 By this we mean that for each point $p \in S$ and $1$-form $\tau \in (T_{d\pi(p)}B)^*$, the projection of the $\Om$-symplectic gradient of $(d \pi |_p)^*(\tau)$ is equal to the $\si$-symplectic gradient of $\tau$. This is the case for instance if 
$$
(M, \om) \hookrightarrow S \stackrel{\pi}{\to} B
$$
is a fibration  with structure group equal to $\Diff^0_{\om}(M, \om)$, the identity component of the symplctic group, such that the $\om_b$-forms on each fiber admit an extension to a closed $2$-form $\tau$ on the total space $S$ that  has constant
rank equal to $\dim M$. Let $\si$ be a
symplectic form on $B$ and take on $S$ the symplectic form
$\Om = \pi^*(\si) + \tau$ (it is indeed symplectic because
the kernel of $\tau$ has constant rank and must be transversal to the
fibers). 

Suppose now that there is a quasi-autonomous Hamiltonian
$H_{t \in [0, T]}:B \to
\R$ whose flow has no non-constant contractible closed orbit  at time $T$ and for which
there are fixed maximum and minimum where the linearization of the flow 
satisfies the conditions of Theorem~\ref{main1}. Then the pull-back of $H_{t  \in [0, T]}$
induces a length minimizing path on $S$. This is obvious since, with our
definition of
$\Om$, the symplectic gradient of $H_t \circ \pi$ projects to the 
symplectic
gradient of $H_t$, thus any non-constant closed $T$-orbit projects to a
non-contractible closed orbit in $B$ and must therefore be 
non-contractible
too. Apply to the two fibers that represent the maximum
and minimum of $H_t \circ \pi$ the same argument on the linearizations and approximate the Hamiltonian near the two fibers so that it become generically under-twisted.
Theorem~\ref{main1} then yields the desired conclusion.

\section{Appendix: Proof of Claim \ref{claim:extra}}\label{sec:appendix}

The following preliminary result implies that ``nearby'' Hamiltonians
 may be used when one considers the length minimizing properties of 
a Hamiltonian path. We refer to Oh \cite{oh3}, Lemma 5.1, for the 
natural optimal $C^0$-version stated here:

\begin{Lemma}\label{near}
Let $\{G_i\} \subset \Hh$ be a sequence of normalized Hamiltonians 
such that each $G_i$ 
generates a path that minimizes the positive Hofer length. If
there is a normalized function $G$ such that the following convergence 
conditions are
satisfied
\begin{enumerate}
\item $\phi^1_{G_i} \to \phi^1_G$ in the $C^0$-norm
\item $ \| G_i - G \| \to 0,$
\end{enumerate}
then $G$ also minimizes the positive Hofer length.
\end{Lemma}

\begin{Remark}
Similar statements apply to the negative and standard Hofer lengths.
\end{Remark}
\begin{Remark}
The convergence criteria of Lemma \ref{near} are satisfied if $G_i \to G$
in the $C^2$-topology.
\end{Remark}

Now, let $H,K \in \Hh$ satisfy the hypotheses of Theorem \ref{main-general}. 
Recall that this means that
\begin{enumerate}
\item $H$ has an under-twisted fixed global maximum at $P$, 
\item $H$ dominates $K$ at $P$,
\item $K$ has a generically under-twisted fixed global maximum at $P$,
\item the flow of $K$ has no contractible $1$-periodic orbits with action 
greater than $\Aa_K(P)$.
\end{enumerate}

To prove Claim \ref{claim:extra}, we must find a pair of functions $H'$ 
and $K'$ such that:
 $H'$ and $K'$ satisfy the hypotheses of Theorem \ref{main-general}, $H'$ 
and $K'$ have the 
additional properties described in Claim \ref{claim:extra}, and $H'$ is 
arbitrarily close to $H$ 
in the sense of Lemma \ref{near}.

To do this, we will perturb $H$ and $K$ several times until the resulting 
functions $H'$ and $K'$ have the properties 
described in Claim \ref{claim:extra}. Each perturbation will be small in 
the sense of Lemma \ref{near} and preserve the 
conditions required by the hypotheses of Theorem \ref{main-general}. The 
only difficulty  occurs in the 
verification that the flow of $K'$ still has no contractible $1$-periodic 
orbits with action greater than $\Aa_{K'}(P)$.

\medskip

  We begin by replacing $H$ and $K$ with the periodic Hamiltonians 
$$
H_{\alpha} = \dot{\alpha}(t)H(\alpha(t),x)\,\,\, \text{    and      
}\,\,\, K_{\alpha} = \dot{\alpha}(t)K(\alpha(t),x)
$$
where $\alpha \colon [0,1] \to [0,1]$ is a smooth function such that 
$\alpha(t)=0$ near $t=0$, $\alpha(t)=1$ near $t=1$,
and $\alpha(t)-t$ is $C^0$-small. These periodic Hamiltonians are 
arbitrarily close to $H$ and $K$ in the 
sense of Lemma \ref{near} and they still satisfy the hypotheses of Theorem 
\ref{main-general}. In fact, the resulting flows are just 
reparameterizations of the original ones 
by $t \to \tau(t)$. Thus, for simplicity, we just may assume that $H$ and 
$K$ are periodic to begin with.  

Now, let $U_P$ be a Darboux chart around $P$. Because $P$ is a generically 
under-twisted fixed
maximum of $K$, it is an isolated fixed point and we can choose $U_P$ to 
be small enough so that none of 
the other $1$-periodic orbits of $\phi^t_K$ enter $U_P$.  We denote the 
canonical norm in this chart by 
$| \,\,\,|$ and the ball of radius $\rho$ by $B_{\rho}(P)$. 

Consider the small bump function 
$$
g(|(\phi^t_{K})^{-1}(x)|^2),
$$
where $g\colon \R \to \R$ is a smooth nonnegative function which strictly 
decreases on 
$[0,\rho^2]$ and vanishes on $[\rho^2, \infty]$. We choose
$\rho$ to be sufficiently small so that $(\phi^t_K)^{-1}(x) \in U_P$ for 
all
$t \in [0,1]$ and $x \in B_{\rho}(P)$.

Set 
$$
H'= H + g(|(\phi^t_K)^{-1}(x)|^2)\,\,\, \text{    and      } \,\,\,K' = K 
+ g(|(\phi^t_K)^{-1}(x)|^2).
$$
When $g$ is sufficiently $C^1$-small, $P$ is still an under-twisted fixed 
global extrema for these new 
Hamiltonians which still clearly satisfy the first three hypotheses of 
Theorem \ref{main-general} and 
now have the second property of Claim \ref{claim:extra}. 

We must now check that $K'$ has no $1$-periodic orbits with action greater 
than $\Aa_{K'}(P)$.
First we prove that the addition of the bump does not create any new 
$1$-periodic orbits.
To do this, we adapt an argument of Siburg's in \cite{si}.

The flow of $K'$ is equal to the composed flow
$\phi^t_{K} \circ \phi^t_g(x)$. The inside flow is $ \phi^t_g(x) 
=
R(t,|x|^2)$, where $R$ is the rotation matrix
$$
R(t,s)=e^{2 g'(s)it} \colon \R^{2n}
\to \R^{2n}.
$$
Note that $\phi^t_{K} \circ \phi^t_g(x)=
\phi^t_{K}(x)$ whenever $x \notin B_{\rho}(P)$. Hence, we only 
need to
prove that there are no new periodic orbits starting at $x \in 
B_{\rho}(P)$.

Let $S_r$ be a sphere of any radius $r$. We define on each $S_r$ the
normalized
distance function, $d_r \colon S_r \times S_r \to \R^+$, which identifies
$S_r$ with
the unit sphere and measures the distance between points with respect to
the usual
round metric.

Consider the lower semi-continuous function
$\mathcal{D}_K \colon B_{\rho}(P) \to \R^+$ given by
$$
\mathcal{D}_K(x) =
  \begin{cases}
    d_{|x|}(x, \phi^1_{K}(x)) & \text{if }  | 
\phi^1_{K}(x)|
= |x|, \\
    1 & \text{otherwise.   }
\end{cases}       
$$

Using the fact that the flow $\phi^t_{K}$ is close to
its generically under-twisted linearized flow near $P$ , it is 
straightforward to 
check that $\mathcal{D}_K$ is bounded away from zero on  $B_{\rho}(P)$, 
for $\rho$ sufficiently small.
The flow $\phi^t_g$ preserves each sphere $S_r$ for $r \in [0,\rho]$. For 
a sufficiently $C^1$-small $g$, we 
can conclude that the perturbed function $\Dd_{K'}$ is also bounded 
away from zero. Hence, the flow $\phi^t_{K'}= \phi^t_{K} \circ 
\phi^t_g(x)$ has no new $1$-periodic orbits.

By our construction, the $1$-periodic orbits of $K'$ are identical to 
those of $K$. 
Furthermore, only the $\Aa_{K'}$-value of $P$ is different than its 
$\Aa_{K}$-value and
$\Aa_{K'}(P) > \Aa_{K}(P)$. So, not only is $\Aa_{K'}(P)$ still a maximum 
in $\Crit (\Aa_{K'})$,
but it is now an {\bf isolated} maximum. 


Finally, we recall that the space $\Hh_{reg}$ is of second category in 
$C^{\infty}(S^1 \times M, \R)$.
We may then perturb $H'$ and $K'$, {\bf away from $P$}, so that the 
resulting functions 
(which we still call $H'$ and $K'$) are in $\Hh_{reg}$. After this final 
perturbation, 
$\Aa_{K'}(P)$ is still an isolated maximum in $\Crit(\Aa_{K'})$. The other 
properties also persist.

\end{document}